\documentclass[11pt,english]{amsart}
\usepackage[latin1]{inputenc}
\usepackage{geometry}
\usepackage{setspace}
\usepackage{amssymb}

\makeatletter

\providecommand{\LyX}{L\kern-.1667em\lower.25em\hbox{Y}\kern-.125emX\@}

 \theoremstyle{plain}    
 \newtheorem{thm}{Theorem}[section]
 \numberwithin{equation}{section} 
 \numberwithin{figure}{section} 
 \theoremstyle{definition}
 \newtheorem{defn}[thm]{Definition}
 \theoremstyle{definition}
  \newtheorem*{example*}{Example}
 \theoremstyle{plain}    
 \newtheorem{lem}[thm]{Lemma} 
 \theoremstyle{plain}    
 \newtheorem*{thm*}{Theorem} 

\usepackage{babel}
\makeatother
\begin{document}

\title{averaging sequences and abelian rank in amenable groups}

\author{Michael Hochman}

\curraddr{Einstein Institute of Mathematics, Edmond J. Safra Campus, Givat
Ram, The Hebrew University of Jerusalem, Jerusalem 91904, Israel }

\email{mhochman@math.huji.ac.il}

\keywords{Amenable group, F\o{}lner sequence, Tempel'man sequence, tempered
sequence}

\begin{abstract}
We investigate the connection between the abelian rank of a countable
amenable group and the existence of good averaging sequences (eg for
the ergodic theorem). 

We show that if $G$ is a group with finite abelian rank $r(G)$,
then $2^{r(G)}$ is a lower bound on the constant associated to a
Tempel'man sequence, and if $G$ is abelain there is a Tempel'man
sequence in $G$ with this constant. On the other hand, infinite rank
precludes the existence of Tempel'man sequences and forces all tempered
sequences to grow super-exponentially.
\end{abstract}
\maketitle

\section{Introduction}

\pagestyle{myheadings}

\markboth{Michael Hochman}{Averaging sequences and abelian rank in
amenable groups}

A countable group $G$ is \emph{amenable} if there exists a sequence
$\left\{ F_{n}\right\} $ of finite subsets of $G$ with the property
that for every $g\in G$,\[
\frac{1}{|F_{n}|}|gF_{n}\cap F_{n}|\rightarrow 1\]
 (See \cite{Green69} for some equivalent definitions). Such a sequence
is called a F\o{}lner sequence. For example, $F_{n}=\left\{ 0,\ldots ,n\right\} \subseteq \mathbb{Z}$
satisfy this condition.

Much of ergodic theory - ie the study of the dynamics of measure-preserving
transformations - can be extended to measure-preserving actions of
amenable groups. A good example of this, which explains why F\o{}lner
sequences are sometimes called averaging sequences, is the mean ergodic
theorem: given a measure preserving action $(g,\omega )\mapsto g\omega $
of $G$ on a probability space $(\Omega ,\mathcal{F},\mu )$, for
any $\varphi \in L^{2}(\Omega ,\mathcal{F},\mu )$ and any F\o{}lner
sequence $\left\{ F_{n}\right\} $ ,\begin{equation}
\frac{1}{|F_{n}|}\sum _{g\in F_{n}}\varphi ^{g}\rightarrow \mathbb{E}(\varphi |\mathcal{F}_{0})\label{eq:mean-ergodic-theorem}\end{equation}
 in $L^{2}$. Here $\mathcal{F}_{0}\subseteq \mathcal{F}$ is the
$\sigma $-algebra of $G$-invariant measurable sets, and $\varphi ^{g}(x)=\varphi (gx)$.
F. Riesz's proof of the mean ergodic theorem for $\mathbb{Z}$-actions
carries over almost verbatim to the general case.

However, for some results - in particular for many pointwise results
- it is necessary to place some restrictions on the averaging sequence
$\left\{ F_{n}\right\} $. One such condition is the following, which
requires that the sets $F_{n}$ be invariant to some degree to translation
by their own elements. As usual, we write $AB=\left\{ ab\, :\, a\in A\, ,\, b\in B\right\} $
and $A^{-1}=\{a^{-1}\, :\, a\in A\}$.

\begin{defn}
\label{def:Templeman-sequence}An increasing F\o{}lner sequence $\left\{ F_{n}\right\} $
in $G$ is called a \emph{Tempel'man sequence} if there is a constant
$C$ such that for every $n$, $|F_{n}^{-1}F_{n}|\leq C|F_{n}|$
\end{defn}
Using Tempel'man sequences in the role of averaging sequence one can
prove the pointwise ergodic theorem (Birkhoff's theorem), which states
that the limit (\ref{eq:mean-ergodic-theorem}) converges almost surely
for $\varphi \in L^{1}$ \cite{Temp67}, and the Shannon-McMillan-Breimann
theorem \cite{OW83}.

Two classes of groups which possess Tempel'man sequences are the locally
finite groups and the finitely generated groups with polynomial growth.
For a long time no examples were known of groups without Tempel'man
sequences. Recently E. Lindenstrauss \cite{Lin01} showed that there
do exist such groups, demonstrating that the lamplighter group $L$
is one such group. This group is defined as follows: let $\mathbb{Z}$
act on $V=\oplus _{i=-\infty }^{\infty }(\mathbb{Z}/2\mathbb{Z})$
by coordinate shift, and set $L=\mathbb{Z}\ltimes V$. 

In this paper we investigate the connection between the abelian rank
of a group and the behavior of the averaging sequences it contains.
Rank is defined as follows:

\begin{defn}
\label{def:rank}The \emph{abelian rank} of a group $G$ is \[
r(G)=\sup \left\{ n\, :\, G\textrm{ contains a subgroup isomorphic to }\mathbb{Z}^{n}\right\} \]

\end{defn}
Our first result connects abelian rank with Tempel'man sequences:

\begin{thm*}
Let $G$ be a countable group. 
\begin{enumerate}
\item If $r(G)<\infty $ then $2^{r(G)}$ is a lower bound on the constant
associated with a Tempel'man sequence in $G$. If in addition $G$
is abelian then there exists in $G$ a Tempel'man sequence with constant
$2^{r(G)}$.
\item If $r(G)=\infty $ then $G$ has no Tempel'man sequences.
\end{enumerate}
\end{thm*}
This provides many new examples of groups without Tempel'man sequences,
for instance $\oplus _{i=-\infty }^{\infty }Z$.

With regard to the problem of the existence of good averaging sequences,
Lindenstrauss demonstrated in \cite{Lin01} that tempered sequences,
defined below, can serve as averaging sequences for the pointwise
ergodic theorem and the SMB theorem:

\begin{defn}
\label{def:tempered}A F\o{}lner sequence $\left\{ F_{n}\right\} $
is \emph{tempered} if there exists a constant $C$ such that $|(\cup _{i<n}F_{i}^{-1})F_{n}|\leq C|F_{n}|$
for every $n$. 
\end{defn}
Note that every Tempel'man sequences is tempered. 

Tempered sequences have the advantage that they exists in every amenable
group; in fact every F\o{}lner sequence has a tempered subsequence.
However tempered sequences may grow quickly. For example, every tempered
sequence in the lamplighter group grows super-exponentially \cite{Lin01}.
Fast growth implies, for instance, that in order to estimate the mean
of a process from observations the size of the data set needed can
grow sharply from one estimate to the next. In contrast, when calculating
the mean of a $\mathbb{Z}$-process using a running average, each
additional symbol of the output can be used to update the estimate.

Our second result shows that abelian rank is relevant to this question
as well:

\begin{thm*}
Let $G$ be a countable group. If $r(G)=\infty $ then $G$ does not
possess Tempel'man sequences, and every tempered sequence in $G$
grows super-exponentially.
\end{thm*}
These results leave open the class of non-abelian finite rank groups,
which contains groups both with and without {}``good'' averaging
sequences. We will say a little about this at the end of section \ref{sec:tempered-sequences}.

Our main tool is the observation that the Brunn-Minkowsky inequality,
which bounds the volumes of sums of sets in $\mathbb{R}^{d}$, can
be used to bound the sums of sufficiently invariant sets in $\mathbb{Z}^{d}$.
This is developed in the next section. In section \ref{sec:tempelman-sequences}
we derive the results about Tempel'man sequences. Section \ref{sec:tempered-sequences}
addresses the growth of tempered sequences.

\medskip{}
\noindent \emph{Acknowledgment.} This work is part of the author's
M.A. thesis, conducted under the guidance of Professor Benjamin Weiss,
whom I would like to thank for all his help and encouragement. I would
also like to thank E. Glasner for raising some of the questions addressed
here.

\section{\label{sec:main-lemma}Tempel'man sequences in integer lattices}

We will need the following classical theorem (see eg \cite{Gardner02}):

\begin{thm*}
\emph{(Brunn-Minkowsky inequality)} Let $A,B\subseteq \mathbb{R}^{d}$
be measurable sets such that $A+B$ is also measurable. Then, denoting
the $d$-dimensional Lebesgue measure by $\nu $,\[
\nu (A+B)\geq \left(\nu (A)^{1/d}+\nu (B)^{1/d}\right)^{d}\]

\end{thm*}
We will use this to obtain a similar statement for sets in $\mathbb{Z}^{d}$
when one of the sets is sufficiently self-invariant, with counting
measure replacing Lebesgue measure. Note that in general the re-statement
of the Brunn-Minkowsky for $\mathbb{Z}^{d}$ is false. This is demonstrated
by considering the {}``simplex'' $A=B=\left\{ 0,e_{1},\ldots ,e_{d}\right\} $
with $e_{i}$ the standard generators of $\mathbb{Z}^{d}$. 

We note that a variant of the Brunn-Minkowsky inequality for integer
lattices was obtained by R.J Gardner and P. Gronchi in \cite{GG01}.
Also note that that I.Z. Ruzsa \cite{Ruzsa94} has obtained inequalities
for sums of sets in $\mathbb{Z}^{d}$ which are sufficient for the
qualitative results stated below. However, the tight bound below does
not follow from that work.

Denote by $e_{i}\in \mathbb{Z}^{d}$ the standard generators of $\mathbb{Z}^{d}$,
that is, $(e_{i})_{j}=\delta _{i,j}$.

\begin{thm}
\label{thm:discrete-brunn-minkowsky}Let $A,B\subseteq \mathbb{Z}^{d}$
and $\delta >0$ and suppose that $A$ is $(1-\delta )$-invariant
to $e_{1},\ldots ,e_{d}$, ie \[
|A\cap (A+e_{i})|>(1-\delta )|A|\]
 for all $i=1,\ldots ,d$. Then \[
|A+B|\geq (1-2k^{2}\delta )\left(|A|^{1/d}+|B|^{1/d}\right)^{d}\]

\end{thm}
\begin{proof}
Let \begin{eqnarray*}
A_{0} & = & \left\{ u\in A\, :\, u+e_{i}\in A\textrm{ and }u+e_{i}+e_{j}\in A\, \textrm{for all }i,j=1,\ldots ,d\right\} \\
 & = & A\; \cap \; \bigcap _{1\leq i\leq d}(A-e_{i})\; \cap \; \bigcap _{1\leq i<j\leq d}(A-e_{j}-e_{i})
\end{eqnarray*}
 From the invariance assumption, $A\cap (A-e_{i})$ and $A\cap (A-e_{i}-e_{j})$
are of size at least $(1-2\delta )|A|$, so\begin{equation}
|A_{0}|>(1-2d^{2}\delta )|A|\label{eq:A-0-size-bound}\end{equation}
 Consider now the sets $\widetilde{A}_{0},\widetilde{B}\subseteq \mathbb{R}^{d}$
obtained by ''thickening'' $A_{0},B$:\[
\widetilde{A}_{0}=A_{0}+[0,1)^{d}\; ,\; \widetilde{B}=B+[0,1)^{d}\]
 (we identify $A_{0},B$ with subsets of $\mathbb{R}^{d}$ in the
obvious way). Clearly \[
\nu (\widetilde{A}_{0})=|A_{0}|\; ,\; \nu (\widetilde{B})=|B|\]
 so by the Brunn-Minkowsky inequality and equation \ref{eq:A-0-size-bound},\begin{eqnarray*}
\nu (\widetilde{A}_{0}+\widetilde{B}) & \geq  & \left(|A_{0}|^{1/d}+|B|^{1/d}\right)^{d}\\
 & \geq  & (1-2d^{2}\delta )\left(|A|^{1/d}+|B|^{1/d}\right)^{d}
\end{eqnarray*}
 On the other hand, it is a simple consequence of the definition of
$A_{0}$ that if $u\in \mathbb{Z}^{d}\cap (\widetilde{A}_{0}+\widetilde{B})$
then $u\in A+B$, and therefore \[
\widetilde{A}_{0}+\widetilde{B}\subseteq (A+B)+[0,1)^{d}\]
 so\[
\nu (\widetilde{A}_{0}+\widetilde{B})\leq |A+B|\]
 putting this all together, we get \[
|A+B|\geq (1-2d^{2}\delta )\left(|A|^{1/d}+|B|^{1/d}\right)^{d}\qedhere \]

\end{proof}
One corollary of this is that in $\mathbb{Z}^{d}$, every Tempel'man
sequence must have constant at least $2^{d}$.

\section{\label{sec:tempelman-sequences}Tempel'man sequences and abelian
rank}

Suppose $G$ is a countable group and contains elements $e_{1},\ldots ,e_{d}$
which generate a subgroup isomorphic to $\mathbb{Z}^{d}$. Then theorem
\ref{thm:discrete-brunn-minkowsky} can be extended to products of
sufficiently invariant sets in $G$, where invariance is measured
relative to $e_{1},\ldots ,e_{d}$.

We first assume that one of the sets is contained in this abelian
subgroup:

\begin{lem}
\label{lem:general-products}Let $G$ be a group and $H\cong \mathbb{Z}^{d}<G$.
Let $A\subseteq H$ be $(1-\delta )$-invariant with respect to the
standard generators $e_{1},\ldots ,e_{d}\in H$, ie $|(e_{i}+A)\cap A|\geq (1-\delta )|A|$
for $i=1,\ldots ,d$. Then for any finite $B\subseteq G$, \[
|BA|,|AB|\geq (1-2d^{2}\delta )\left(|A|^{1/d}+|B|^{1/d}\right)^{d}\]

\end{lem}
\begin{proof}
We show the inequality for $|AB|$. Let $\left\{ Hg_{i}\right\} $
be an enumeration of the right cosets of $H$ and set $B_{i}=B\cap Hg_{i}$.
Using theorem \ref{thm:discrete-brunn-minkowsky} ,\begin{eqnarray*}
|AB| & = & \sum _{i}|AB_{i}|\\
 &  & \geq (1-2d^{2}\delta )\sum _{i}(|A|^{1/d}+|B_{i}|^{1/d})^{d}
\end{eqnarray*}
 using the fact that $t\mapsto (c+t^{1/d})^{d}$ is concave and homogeneous
and that $\sum _{i}|B_{i}|=|B|$,\[
\geq (1-2d^{2}\delta )(|A|^{1/d}+|B|^{1/d})^{d}\qedhere \]

\end{proof}
Next we drop the requirement that the invariant set be contained in
the abelian subgroup, and instead consider products of the form $F^{-1}F$:

\begin{lem}
\label{lem:same-size-products}Suppose $G$ is a group, $H\cong \mathbb{Z}^{d}<G$
with $e_{1},\ldots ,e_{d}$ the standard generators in $H$, and $F\subseteq G$
finite such that $F$ is $(1-\delta )$-invariant with respect to
left multiplication by $e_{1},\ldots ,e_{d}$, ie $\frac{1}{|F|}|e_{i}F\cap F|\geq (1-\delta )|F|$.
Then \[
|F^{-1}F|\geq 2^{d}(1-2d^{2}\sqrt{\delta })(1-d\sqrt{\delta })|F|\]

\end{lem}
\begin{proof}
Let $Hg_{i}$ be an enumeration of the right cosets of $H$ and $F_{i}=F\cap Hg_{i}$.
From the invariance assumption there is a set of indices $I$ such
that \[
|\cup _{i\in I}F_{i}|\geq (1-d\sqrt{\delta })|F|\]
 and $F_{i}$ is $(1-\sqrt{\delta })$-invariant to $e_{1},\ldots ,e_{d}$
from the left for each $i\in I$. Let $i_{0}\in I$ be such that $|F_{i_{0}}|$
is maximal among $\left\{ |F_{i}|\right\} _{i\in I}$. Now\[
|F^{-1}F|\geq |F_{i_{0}}^{-1}F|\geq \sum _{i\in I}|F_{i_{0}}^{-1}F_{i}|\geq \]
 by the lemma \ref{lem:general-products},\[
\geq \sum _{i\in I}(1-2d^{2}\sqrt{\delta })(|F_{i_{0}}|^{1/d}+|F_{i}|^{1/d})^{d}\geq \]
 by choice of $i_{0}$\[
\geq (1-2d^{2}\sqrt{\delta })\sum _{i\in I}2^{d}|F_{i}|\geq 2^{d}(1-2d^{2}\sqrt{\delta })(1-d\sqrt{\delta })|F|\qedhere \]

\end{proof}
Recall that the abelian rank $r(G)$ of $G$ is the largest $n$ such
that $\mathbb{Z}^{d}$ can be embedded in $G$ (definition \ref{def:rank}).
Recall also the definition of Tempel'man sequences (definition \ref{def:Templeman-sequence}).
As an immediate corollary of lemma \ref{lem:same-size-products} we
have

\begin{thm}
If $r(G)<\infty $ then $2^{r(G)}$ is a lower bound on the constant
associated with Tempel'man sequences in $G$.
\end{thm}
For abelian groups, there is a converse to this, which together with
theorem \ref{thm:no-Templeman-sequences} characterizes those abelian
groups which have Tempel'man sequences. I thank B. Weiss for suggesting
the following construction:

\begin{thm}
\label{thm:finite-rank-abelian}If $A$ is a countable abelian group
with $r(A)<\infty $ then there are Tempel'man sequences in $A$. 
\end{thm}
\begin{proof}
Suppose that the rank of $A$ is $d<\infty $. We will construct a
Tempel'man sequence. Let $\left\{ a_{1},a_{2},\ldots \right\} $ be
an enumeration of the elements of $A$. Let $A_{n}$ be the group
generated by $\left\{ a_{1},\ldots ,a_{n}\right\} $. Since $A_{n}$
is finitely generated, we can write it as $T_{n}\oplus \mathbb{Z}^{d(n)}$
with $T_{n}$ a finite group and $d(n)\leq d$. Set $F_{n}=T_{n}\times \left\{ 0,\ldots ,k(n)\right\} ^{d(n)}$,
where $k(n)$ is chosen large enough that $F_{n}$ is $\frac{1}{n}$-invariant
to $a_{1},\ldots ,a_{n}$ and $F_{n-1}\subseteq F_{n}$ (it is easy
to check that such a $k(n)$ exists). Now one verifies that $|F_{n}-F_{n}|\leq 2^{d(n)}|F_{n}|$.
Thus $\left\{ F_{n}\right\} $ is a Tempel'man sequence. 
\end{proof}
From the proof one also sees that if $r(A)=d<\infty $ then there
exist Tempel'man sequences with constant $2^{d}$, and lemma \ref{lem:same-size-products}
shows that is the best possible.

On the other hand, as noted in the introduction, finite rank in general
does not guarantee the existence of Tempel'man sequences: the lamplighter
group $L$ studied in \cite{Lin01} has $r(L)=1$ but has no Tempel'man
sequences.

From lemma \ref{lem:same-size-products} we also derive

\begin{thm}
\label{thm:no-Templeman-sequences}If $G$ is an amenable group with
$r(G)=\infty $ then $G$ does not have Tempel'man sequences. 
\end{thm}

\section{\label{sec:tempered-sequences}Growth of tempered sequences}

The previous section and the work of Lindenstrauss on the lamplighter
group \cite{Lin01} show that there are amenable groups without Tempel'man
sequences. However as noted in the introduction all amenable groups
possess tempered sequences (definition \ref{def:tempered}), which
can serve as averaging sequences. Ideally, we would like there to
be averaging sequences which grow slowly (in $\mathbb{Z}$ they can
grow linearly). It is therefore of interest to find conditions indicating
or precluding such slow growth.

In order to analyze the growth rate of tempered sequences we will
need information about the product of different sets, each of which
is somewhat invariant:

\begin{lem}
\label{lem:different-size-products}Suppose $G$ is a group, $H\cong \mathbb{Z}^{d}<G$
with $e_{1},\ldots ,e_{d}$ the standard generators in $H$, and $F_{1,}F_{2}\subseteq G$
finite such that $F_{1},F_{2}$ are $(1-\delta )$-invariant with
respect to left-multiplication by $e_{1},\ldots ,e_{d}$, ie $|e_{i}F_{k}\cap F_{k}|\geq (1-\delta )|F_{k}|$
for $i=1,\ldots ,d$ and $k=1,2$. Then \[
|F_{1}^{-1}F_{2}|\geq 2^{d}(1-2d^{2}\sqrt{\delta })(1-d\sqrt{\delta })\min \left\{ |F_{1}|,|F_{2}|\right\} \]

\end{lem}
\begin{proof}
Similar to lemma \ref{lem:same-size-products}. Let $Hg_{i}$ be an
enumeration of the right cosets of $H$ and $F_{1,i}=F_{1}\cap Hg_{i}$
and similarly $F_{2,i}$. From the invariance assumption there are
sets of indices $I,J$ such that \[
|\cup _{i\in I}F_{1,i}|\geq (1-d\sqrt{\delta })|F_{1}|\; ,\; |\cup _{j\in J}F_{2,j}|\geq (1-d\sqrt{\delta })|F_{2}|\]
 and $F_{k,i}$ is $(1-\sqrt{\delta })$-invariant with respect to
left multiplication by $e_{1},\ldots ,e_{d}$ for each $i\in I$ and
$k=1,2$. Let $i_{0}\in I$ be such that $|F_{1,i_{0}}|$ is maximal
among $\left\{ |F_{1,i}|\right\} _{i\in I}$ and $j_{0}\in J$ such
that $|F_{2,j_{0}}|$ is maximal among $\left\{ |F_{2,j}|\right\} _{j\in J}$.
Now the calculation in lemma \ref{lem:same-size-products} shows that
if $|F_{1,i_{0}}|\geq |F_{2,j_{0}}|$ then \[
|F_{1}^{-1}F_{2}|\geq 2^{d}(1-2d^{2}\sqrt{\delta })(1-d\sqrt{\delta })|F_{2}|\]
 while if $|F_{2,j_{0}}|\geq |F_{1,i_{0}}|$ then\[
|F_{1}^{-1}F_{2}|\geq 2^{d}(1-2d^{2}\sqrt{\delta })(1-d\sqrt{\delta })|F_{1}|\]
 which together complete the lemma. 
\end{proof}
\begin{thm}
\label{thm:temepered-sequences}If $G$ is amenable and $r(G)=\infty $
then every tempered sequence in $G$ grows super-exponentially. 
\end{thm}
\begin{proof}
Let $\left\{ F_{n}\right\} $ be a tempered F\o{}lner sequence in
$G$, so for some constant $C$, $|F_{n-1}^{-1}F_{n}|\leq C|F_{n}|$
for all $n$. Fix $d>2+\log _{2}C$ and $\delta =\frac{1}{16d^{2}}$.
Let $H\cong \mathbb{Z}^{d}<G$ be generated by $e_{1},\ldots ,e_{d}$.
For large enough $n$ the $F_{n}$'s will be $(1-\delta )$-invariant
to $e_{1},\ldots ,e_{d}$, so by lemma \ref{lem:different-size-products}
for large enough $n$ we will have \[
C|F_{n}|\geq |F_{n-1}^{-1}F_{n}|\geq 2^{d}(1-2d^{2}\sqrt{\delta })(1-d\sqrt{\delta })\min \left\{ |F_{n}|,|F_{n-1}|\right\} \]
 substituting $\delta =\frac{1}{16d^{2}}$ and rearranging gives\[
|F_{n}|\geq \frac{2^{d-2}}{C}\min \{|F_{n}|,|F_{n-1}|\}\]
 Since $d>2+\log _{2}C$, this relation excludes the possibility that
$\min \{|F_{n-1}|,|F_{n}|\}=|F_{n}|$, so it must be that \[
|F_{n}|\geq \frac{2^{d-2}}{C}|F_{n-1}|\]
 Since $d$ was arbitrary and this holds for all large enough $n$,
the proof is complete. 
\end{proof}
The results presented in this paper do not settle the question of
the character of averaging sequences in the class nonabelian groups
of finite rank. In this class we have examples of groups with good
sequences (eg nonabelian groups with polynomial growth) and groups
with bad sequences (eg the lamplighter group $L$). Lindenstrauss
has speculated that for finitely generated groups exponential growth
may force super-exponential growth of tempered sequences the absence
of Tempel'man sequences. The following example provides some further
support for this. Note that unlike $L$, it is torsion-free.

\begin{example*}
Let $\mathbb{Z}$ act on $\oplus _{i=-\infty }^{\infty }\mathbb{Z}$
by shift to the left, ie for $n\in \mathbb{Z}$ and $u=\left(u_{i}\right)\in \oplus _{i=-\infty }^{\infty }\mathbb{Z}$
we define $\left(u^{n}\right)_{i}=u_{i+n}$. Let \[
G=\mathbb{Z}\ltimes \oplus _{i=-\infty }^{\infty }\mathbb{Z}\]
 explicitly, $G=\left\{ (n,u)\in \mathbb{Z}\times \oplus _{i=-\infty }^{\infty }\mathbb{Z}\right\} $
and $(n,u)\cdot (m,v)=(n+m,u^{m}+v)$. 

This group is finitely generated (eg. by $(0,(\ldots ,0,0,1,0,0,\ldots ))$
and $(1,(\ldots ,0,0,0,\ldots )$) and is torsion-free. $G$ is also
amenable, since it is solvable \cite{Green69}. Clearly $r(G)=\infty $,
so by theorem \ref{thm:temepered-sequences} $G$ has no Tempel'man
sequences and every tempered sequence in $G$ grows super-exponentially.
\end{example*}
\bibliographystyle{plain}
\bibliography{teza}

\end{document}